\documentclass[12pt]{amsart}
\usepackage{a4}
\usepackage{amssymb}

\newcommand{\R}{{\mathbb R}}

\newcommand{\Z}{{\mathbb Z}}

\newcommand{\LA}{{\mathfrak a}}

\newcommand{\LG}{{\mathfrak g}}

\newcommand{\LK}{{\mathfrak k}}

\newcommand{\LM}{{\mathfrak m}}
\newcommand{\LN}{{\mathfrak n}}
\newcommand{\LP}{{\mathfrak p}}
\newcommand{\LQ}{{\mathfrak q}}
\newcommand{\LS}{{\mathfrak s}}

\newcommand{\Sum}{\displaystyle \sum}
\newcommand{\ad}{\mbox{\rm ad}}

\newcommand{\tr}{\mbox{\rm tr}}
\newcommand{\Ric}{\mbox{\rm Ric}}
\newcommand{\ric}{\mathrm{ric}}

\newcommand{\inner}[2]{\langle #1 , #2 \rangle} 

\newtheorem{theorem}{Theorem}[section]
\newtheorem{Prop}[theorem]{Proposition}
\newtheorem{Thm}[theorem]{Theorem}
\newtheorem{Lem}[theorem]{Lemma}

\theoremstyle{definition}

\newtheorem{Def}[theorem]{Definition}

\theoremstyle{remark}

\numberwithin{equation}{section}

\begin{document}

\title[Parabolic subgroups and Einstein solvmanifolds]{Parabolic subgroups of semisimple Lie groups and Einstein solvmanifolds}

\author{Hiroshi Tamaru}
\address{Department of Mathematics, Hiroshima University, 
Higashi-Hiroshima 739-8526, Japan}
\email{tamaru@math.sci.hiroshima-u.ac.jp}
\thanks{The author was supported in part by 
Grant-in-Aid for Young Scientists (B) 17740039, 
The Ministry of Education, Culture, Sports, Science and Technology, Japan.}

\subjclass[2000]{
Primary 53C30; 
Secondary 22E25 
}
\keywords{
Parabolic subgroup, 
Langlands decomposition, 
nilpotent and solvable Lie group, 
Einstein solvmanifold
}

\begin{abstract}
In this paper, we study the solvmanifolds constructed from 
any parabolic subalgebras of any semisimple Lie algebras. 
These solvmanifolds are naturally homogeneous submanifolds 
of symmetric spaces of noncompact type. 
We show that the Ricci curvatures of our solvmanifolds coincide with 
the restrictions of the Ricci curvatures of the ambient symmetric spaces. 
Consequently, all of our solvmanifolds are Einstein, 
which provide a large number of new examples of 
noncompact homogeneous Einstein manifolds. 
We also show that our solvmanifolds are minimal, 
but not totally geodesic submanifolds of symmetric spaces. 
\end{abstract}

\maketitle

\setlength{\baselineskip}{15pt}

\section{Introduction}

A \textit{solvmanifold} is a Riemannian manifold which has 
transitive solvable group of isometries. 
Recently a study of solvmanifolds is active, 
since solvmanifolds provide a lot of interesting examples of Riemannian manifolds 
with some particular geometric properties. 
Most successful examples are Damek-Ricci spaces, 
which provide non-symmetric harmonic manifolds 
(\cite{DR}, see also \cite{BTV})
and hence give counterexamples for Lichnerowicz conjecture. 
Solvmanifolds also provide examples of noncompact homogeneous Einstein manifolds, 
which is the subject of this paper. 
Damek-Ricci spaces and symmetric spaces of noncompact type 
are typical examples of Einstein solvmanifolds. 
It is conjectured that 
every noncompact homogeneous Einstein manifold is a solvmanifold 
(Alekseevskii conjecture, see \cite{Besse}). 
We refer \cite{L-dga} and the references therein 
for known examples of Einstein solvmanifolds. 
The structures of Einstein solvmanifolds have been deeply studied by 
Heber (\cite{H}). 
Further progress has been made recently by Lauret (\cite{L-standard}). 

\medskip 

The purpose of this paper is to construct new examples of Einstein solvmanifolds. 
The starting point of our construction is 
a parabolic subalgebra of a semisimple Lie algebra $\LG$, 
which will be mentioned in Section \ref{section_para}. 
Let $\LQ_{\Lambda'}$ be the parabolic subalgebra of $\LG$ 
corresponding to a subset $\Lambda'$ 
of a set of simple roots $\Lambda$ of the restricted root system of $\LG$. 
Then, $\LQ_{\Lambda'}$ has a decomposition 
$\LQ_{\Lambda'} = \LM_{\Lambda'} + \LA_{\Lambda'} + \LN_{\Lambda'}$ 
into a reductive subalgebra $\LM_{\Lambda'}$ 
and a solvable subalgebra $\LS_{\Lambda '} = \LA_{\Lambda '} + \LN_{\Lambda '}$, 
which is called the Langlands decomposition. 
Our solvmanifolds are simply-connected Lie groups with Lie algebras 
$\LS_{\Lambda '}$, 
equipped with natural left-invariant Riemannian metrics. 
In this paper, we always identify our solvmanifolds 
with the underlying metric solvable Lie algebras 
$(\LS_{\Lambda '} = \LA_{\Lambda '} + \LN_{\Lambda '}, \inner{}{})$. 
Since there is a close relation between parabolic subalgebras and 
semisimple graded Lie algebras (see Section \ref{section_para}), 
our class is a generalization of the Einstein solvmanifolds studied 
by Mori (\cite{Mori}) and the author (\cite{T7}, \cite{T-dga}). 
In fact, Mori (\cite{Mori}) showed that, 
if $\LG$ is a complex simple (classical) Lie algebra, $\dim \LA_{\Lambda '} = 1$, 
and $\LN_{\Lambda '}$ is two-step nilpotent, 
then our solvmanifolds are Einstein. 
The author proved in \cite{T7} that, 
if $\LN_{\Lambda '}$ is two-step nilpotent or $\dim \LA_{\Lambda '} = 1$, 
then they are Einstein. 
The main result of this paper states that, without any assumptions, 
all of our solvmanifolds are Einstein. 

\medskip 

Our class of Einstein solvmanifolds 
$(\LS_{\Lambda '} = \LA_{\Lambda '} + \LN_{\Lambda '}, \inner{}{})$ 
has the following features: 

\medskip 

Firstly, it provides 
a large number of new examples of Einstein solvmanifolds. 
For each semisimple Lie algebra $\LG$, 
parabolic subalgebras are parameterized by subsets $\Lambda '$ of 
the set of simple roots $\Lambda$. 
Denote by $r$ the (split) rank. 
Therefore, one has approximately $2^r-1$ parabolic subalgebras of $\LG$ 
(to be exact, one has to take account of symmetries of the Dynkin diagram). 

\medskip 

Secondly, for our Einstein solvmanifolds 
$(\LS_{\Lambda '} = \LA_{\Lambda '} + \LN_{\Lambda '}, \inner{}{})$, 
the degree of nilpotency of $\LN_{\Lambda '}$ can be arbitrary large. 
One can determine the degree of nilpotency of $\LN_{\Lambda '}$ 
by looking at the coefficients of the highest root. 
In general, it seems to be difficult to find explicit examples 
of Einstein solvmanifolds 
with nilradicals of high nilpotency (cf. \cite{Niko}). 

\medskip 

Thirdly, our class of Einstein solvmanifolds contains 
all symmetric spaces of noncompact type. 
In fact, if we start from the parabolic subalgebra $\LQ_{\emptyset}$ 
corresponding to the empty set $\emptyset \subset \Lambda$, 
then the constructed solvmanifold 
$(\LS_{\emptyset} = \LA_{\emptyset} + \LN_{\emptyset}, \inner{}{})$ 
is isometric to the symmetric space (Proposition \ref{prop:symm}). 
Thus, our class of Einstein solvmanifolds 
is a generalization of symmetric spaces of noncompact type, 
and may be regarded as a higher rank analogue of the study of Damek-Ricci spaces 
(recall that, the class of Damek-Ricci spaces is a generalization 
of noncompact symmetric spaces of rank one). 

\medskip 

Fourthly, 
our solvmanifolds have a remarkable property as Riemannian submanifolds. 
Our solvmanifold 
$(\LS_{\Lambda '} = \LA_{\Lambda '} + \LN_{\Lambda '}, \inner{}{})$ 
is naturally a Riemaniann submanifold of 
the symmetric space of noncompact type 
$(\LS_{\emptyset} = \LA_{\emptyset} + \LN_{\emptyset}, \inner{}{})$. 
In the proof of our main theorem, 
we compare the Ricci curvatures of these manifolds, 
denoted by $\ric^{\LS_{\Lambda '}}$ and $\ric^{\LS_{\emptyset}}$ respectively, 
and show that 
$\ric^{\LS_{\Lambda '}} = 
\ric^{\LS_{\emptyset}} |_{\LS_{\Lambda '} \times \LS_{\Lambda '}}$. 
This immediately yields that all of our solvmanifolds are Einstein, 
since $(\LS_{\emptyset}, \inner{}{})$ is Einstein. 
The above mentioned property, having the same Ricci curvature, seems be interesting. 
A study of submanifolds with this property 
may provide new examples of Einstein solvmanifolds. 
We also show in Section \ref{sec:minimal} that 
$(\LS_{\Lambda '}, \inner{}{})$ are 
minimal, but not totally geodesic submanifolds of 
$(\LS_{\emptyset}, \inner{}{})$. 

\medskip 

The author would like to express his deep gratitude to 
Professor Soji Kaneyuki, 
Professor Toshiyuki Kobayashi, 
and Professor Yusuke Sakane, 
for their useful advice and valuable comments. 

\section{Preliminaries}
\label{section_solv}

In this section, we recall some curvature formulae for 
a simply-connected Lie group $G$ with a left invariant metric $g$. 
A Lie algebra endowed with an inner product is usually called a 
\textit{metric Lie algebra}. 
Curvatures of $(G,g)$ can be completely determined by 
the underlying metric Lie algebra $(\LG, \inner{}{})$. 

\medskip 

We identify the Lie algebra $\LG$ with the set of left-invariant vector fields 
of $G$. 
Let $X, Y \in \LG$. 
The Kozsul formula for Levi-Civita connection $\nabla$ of $(G, g)$ yields that 
\[
\nabla_X Y = (1/2) [X,Y] + U(X,Y) , 
\]
where $U : \LG \times \LG \rightarrow \LG$ denotes the symmetric bilinear form 
defined by 
\[
2 \inner{U(X,Y)}{Z} = \inner{[Z,X]}{Y} + \inner {X}{[Z,Y]} . 
\]
The Riemannian curvature 
$R : \LG \times \LG \times \LG \rightarrow \LG$ 
is defined by 
\[
R(X,Y)Z := \nabla_X \nabla_Y Z - \nabla_Y \nabla_X Z - \nabla_{[X,Y]} Z . 
\]
Let $\{ E_i \}$ be an orthonormal basis of $\LG$. 
The Ricci curvature 
$\ric : \LG \times \LG \rightarrow \R$ is defined by 
\[
\ric(X,Y) := \sum_i \inner{R(X, E_i) Y}{E_i} . 
\]
A metric Lie algebra $(\LG, \inner{}{})$, identified with $(G, g)$, 
is said to be \textit{Einstein} 
if $\ric = c \inner{}{}$ holds for some $c \in \R$. 
If $(\LG, \inner{}{})$ is Einstein, 
the number $c$ is called the \textit{Einstein constant}. 

\medskip 

We have a formula for the Ricci curvature in terms of 
the Killing form $B$ of $\LG$ and the vector $H_0 := \sum U(E_i, E_i) \in \LG$. 

\begin{Prop}[\cite{Besse}, p185, Corollary 7.38]
\label{ricci-formula}
The Ricci curvature $\ric$ of $(\LG, \inner{}{})$ is given by
\begin{eqnarray*}
\ric (X,Y) & = & 
-(1/2) \sum \inner{[X, E_i]}{[Y, E_i]} - (1/2) B(X,Y) \\ 
&& + (1/4) \sum \inner{[E_i,E_j]}{X} \inner{[E_i,E_j]}{Y} - \inner{U(X,Y)}{H_0} . 
\end{eqnarray*}
\end{Prop}

We are interested in solvable Lie algebras, which are called Iwasawa-type. 
Recall that, a Lie algebra $\LS$ is said to be \textit{solvable} 
if the derived subalgebra $[\LS, \LS]$ is nilpotent. 

\begin{Def}
A metric solvable Lie algebra $(\LS, \inner{}{})$ 
is called \textit{standard} if $\LA := [\LS, \LS]^{\perp}$ is abelian. 
We denote by $\LN := [\LS, \LS]$ and thus $\LS = \LA + \LN$. 
A standard solvable Lie algebra 
$(\LS = \LA + \LN, \inner{}{})$ is said to be of \textit{Iwasawa-type} if 
\begin{enumerate}
\item
$\ad_A$ is symmetric for every $A \in \LA$, and 
\item
for some $A_0 \in \LA$, the restriction $\ad_{A_0} |_{\LN}$ is positive definite. 
\end{enumerate}
\end{Def}

For a solvable Lie algebra of Iwasawa-type $(\LS = \LA + \LN, \inner{}{})$, 
the Ricci curvature can be represented in terms of the Ricci curvature 
$\ric^{\LN}$ of $(\LN, \inner{}{})$ and the vector $H_0$ defined above. 
In this case, $H_0$ coincides with the mean curvature vector of the submanifold 
$(\LN, \inner{}{}) \subset (\LS, \inner{}{})$, 
and $H_0 \in \LA$ holds. 
The vector $H_0$ is called the \textit{mean curvature vector} 
of $(\LS = \LA + \LN, \inner{}{})$. 

\begin{Thm}[\cite{W}]
\label{wolter}
Let $(\LS = \LA + \LN, \inner{}{})$ 
be a solvable Lie algebra of Iwasawa-type. 
Then the Ricci curvature $\ric$ satisfies 
\begin{enumerate}
\item 
$\ric(A,A') = - \tr (\ad_A) \circ (\ad_{A'})$ \quad 
for all $A, A' \in \LA$, 
\item 
$\ric(A,X) = 0$ \quad 
for all $A \in \LA$ and $X \in \LN$, 
\item
$\ric(X,Y) = \ric^{\LN}(X,Y) - \inner{\ad_{H_0}X}{Y}$ \quad 
for all $X, Y \in \LN$. 
\end{enumerate}
\end{Thm}

Note that, $\LS' := \R H_0 + \LN$ is a subalgebra of $\LS = \LA + \LN$, 
and $H_0$ is also the mean curvature vector of 
$(\LS', \inner{}{})$. 
Therefore, Theorem \ref{wolter} yields that, 
the Ricci curvature of $(\LS', \inner{}{})$ 
coincides with the restriction of the Ricci curvature of $(\LS, \inner{}{})$. 
Hence, if $(\LS, \inner{}{})$ is Einstein, 
then so is $(\LS', \inner{}{})$. 
The obtained solvmanifold $(\LS', \inner{}{})$ 
is called the \textit{rank one reduction} (see \cite{H}). 

\medskip 

Finally, we consider a metric nilpotent Lie algebra $(\LN, \inner{}{})$. 
In this case, one can see that $B = 0$ and $H_0 = 0$, 
which makes the formula for Ricci curvature simpler. 
We express the formula in terms of the Ricci $(1,1)$-tensor 
$\Ric : \LN \rightarrow \LN$, defined by 
\[
\inner{\Ric(X)}{Y} = \ric(X,Y) . 
\]

\begin{Prop}[\cite{A}]
\label{ricci-nil}
The Ricci tensor $\Ric^{\LN}$ 
of a metric nilpotent Lie algebra $(\LN, \inner{}{})$ is given by 
\[
\Ric^{\LN} = (1/4) \sum \ad_{E_i} \circ (\ad_{E_i})^{\ast} 
- (1/2) \sum (\ad_{E_i})^{\ast} \circ \ad_{E_i} , 
\]
where $\{ E_i \}$ is an orthonormal basis of $\LN$, and 
$(\ad_{E_i})^{\ast}$ denotes the dual of $\ad_{E_i}$ defined by 
$\inner{(\ad_{E_i})^{\ast}(X)}{Y} = \inner{X}{\ad_{E_i}(Y)}$. 
\end{Prop}

\section{Parabolic subalgebras of semisimple Lie algebras}
\label{section_para}

In this section, we review a structure theory of semisimple Lie algebras. 
We mention semisimple graded Lie algebras, parabolic subalgebras, 
and Langlands decompositions. 

\medskip 

Let $\LG$ be a semisimple Lie algebra and $\sigma$ a Cartan involution. 
This gives us a Cartan decomposition 
$\LG = \LK + \LP$, 
where $\LK$ is a maximal compact subalgebra of $\LG$. 
The Cartan involution $\sigma$ also determines a positive definite 
$\ad_{\LK}$-invariant inner product $B_{\sigma}$ on $\LG$ by 
\[
B_{\sigma}(X,Y) := - B(X, \sigma(Y)) , 
\]
where $B$ denotes the Killing form of $\LG$. 
One can easily see that 

\begin{Prop}
\label{prop:inner}
$B_{\sigma}([Z,X],Y) = - B_{\sigma}(X, [\sigma(Z),Y])$ 
for every $X,Y,Z \in \LG$. 
\end{Prop}

Let $\LA$ be a maximal abelian subspace of $\LP$. 
In the usual way, $\LA$ defines the (restricted) root system 
$\Delta = \Delta(\LG, \LA)$ 
of $\LG$ with respect to $\LA$. 
Denote by $\LG_{\alpha}$ the root space of a root $\alpha$. 
We thus obtain the \textit{root space decomposition}, 
\[
\LG = \LG_0 + \sum_{\alpha \in \Delta} \LG_{\alpha} , 
\]
where $\LG_0$ is the centralizer of $\LA$ in $\LG$. 
It is easy to see $\LG_0 = \LK_0 + \LA$, 
where $\LK_0$ is the centralizer of $\LA$ in $\LK$. 

\medskip 

By the inner product $B_{\sigma}$ and the subspace $\LA$, 
one can define root vectors $H_{\alpha}$. 
For $\alpha \in \Delta$, a vector $H_{\alpha} \in \LA$ is called the 
\textit{root vector} of $\alpha$ if 
$B_{\sigma}(H_{\alpha}, A) = \alpha(A)$ holds for every $A \in \LA$. 
For a unit vector $X_{\alpha} \in \LG_{\alpha}$, 
we have $[\sigma(X_{\alpha}), X_{\alpha}] \in \LA$, 
and hence Proposition \ref{prop:inner} yields that  
$H_{\alpha} = [\sigma(X_{\alpha}), X_{\alpha}]$. 

\medskip 

Let $\Lambda$ be a set of simple roots of $\Delta$. 
We call a root $\alpha$ \textit{positive} 
if $\alpha$ can be represented as a linear combination of $\Lambda$ 
with non-negative coefficients. 
Denote by $\Delta^{+}$ the set of positive roots. 
One can easily see that 
$\LN := \sum_{\alpha \in \Delta^+} \LG_{\alpha}$
is a nilpotent subalgebra. 
The decomposition $\LG = \LK + \LA + \LN$ is called an 
\textit{Iwasawa decomposition}. 

\medskip 

Now we define semisimple graded Lie algebras. 
We refer \cite{KA} (and also \cite{T4}). 

\begin{Def}
A decomposition of a semisimple Lie algebra $\LG$ into subspaces, 
\[
\LG = \sum_{k \in \Z} \LG^k , 
\]
is called a \textit{gradation} if 
$[\LG^i, \LG^j] \subset \LG^{i+j}$ holds for every $i,j \in \Z$. 
A graded Lie algebra is said to be of \textit{$\nu$-th kind} 
if $\LG^k = 0$ for $|k| > \nu$. 
\end{Def}

Every gradation can be constructed from a particular element $Z \in \LA$. 
Denote by $\Lambda = \{ \alpha_1, \ldots, \alpha_r \}$, and 
let $\{ H^1, \ldots, H^r \}$ be the dual basis of $\Lambda$, 
that is, $\alpha_i(H^j) = \delta_{ij}$. 

\begin{Prop}[\cite{KA}]
\label{gradation}
Let $Z := c_1 H^1 + \cdots + c_r H^r$ with $c_1, \ldots, c_r \in \Z_{\ge 0}$. 
Then, every eigenvalue of $\ad_Z$ is an integer 
and the eigenspace decomposition $\LG = \sum_{k \in \Z} \LG^k$ gives a gradation. 
Conversely, every gradation can be obtained in this way up to conjugation. 
\end{Prop}

The above element $Z$ is called a \textit{characteristic element}. 
Let $\LG = \sum \LG^k$ be the graded Lie algebra 
with characteristic element $Z \in \LA$. 
By construction, we have 
\[
\LG^0 = \LG_0 + \sum_{\alpha(Z)=0} \LG_{\alpha} , \quad 
\LG^k = \sum_{\alpha(Z)=k} \LG_{\alpha} 
\ \mbox{(for $k \neq 0$)}. 
\]

Our next object is a parabolic subalgebra of a semisimple Lie algebra. 
We refer \cite{Knapp} and \cite{OV}. 

\begin{Def}
The subalgebra $\LK_0 + \LA + \LN$ is called a 
\textit{minimal parabolic subalgebra} of $\LG$. 
A subalgebra $\LQ$ of $\LG$ is called \textit{parabolic} if 
$\LQ$ contains the minimal parabolic subalgebra up to conjugation. 
\end{Def}

Every parabolic subalgebra can be constructed from 
a subset $\Lambda '$ of $\Lambda$. 
Denote by $\langle \Lambda ' \rangle$ the set of the roots spanned by $\Lambda '$.

\begin{Prop}[\cite{Knapp}, Proposition 7.76]
\label{Knapp}
Let $\Lambda'$ be a proper subset of $\Lambda$, and define 
\begin{eqnarray*}
\LQ_{\Lambda '} & := & 
\LG_0 + \sum_{\beta \in \Delta^+ \cup \langle \Lambda ' \rangle} \LG_{\beta} . 
\end{eqnarray*}
Thus, $\LQ_{\Lambda'}$ is a parabolic subalgebra of $\LG$. 
Conversely, every parabolic subalgebra can be constructed in this way 
up to conjugation. 
\end{Prop}

Note that, when $\Lambda' = \emptyset$, the parabolic subalgebra 
$\LQ_{\emptyset}$ is nothing but the minimal parabolic subalgebra. 
There is a correspondence between graded Lie algebras and parabolic subalgebras. 

\begin{Prop}
For a graded Lie algebra $\LG = \sum \LG^k$, 
the subalgebra $\sum_{k \ge 0} \LG^k$ is a parabolic subalgebra of $\LG$. 
Conversely, for every parabolic subalgebra $\LQ$ of $\LG$, 
there exists a gradation $\LG = \sum \LG^k$ such that $\LQ = \sum_{k \ge 0} \LG^k$. 
\end{Prop}

\begin{proof}
Let $\LG = \sum \LG^k$ be a graded Lie algebra. 
Proposition \ref{gradation} states that 
there exists the characteristic element 
$Z = c_{i_1} H^{i_1} + \cdots + c_{i_k} H^{i_k} \in \LA$ 
with $c_{i_1}, \ldots, c_{i_k} \in \Z_{> 0}$, 
up to conjugation. 
Let 
\[
\Lambda' := \{ \alpha_i \in \Lambda \mid \alpha_i(Z) = 0 \} 
= \Lambda \setminus \{ \alpha_{i_1}, \ldots, \alpha_{i_k} \} . 
\]
It is easy to see that 
\begin{eqnarray}
\label{eq:root1}
\langle \Lambda ' \rangle 
& = & \{ \beta \in \Delta \mid \beta(Z) = 0 \} , \\ 
\label{eq:root2}
\Delta^+ \setminus \langle \Lambda ' \rangle 
& = & \{ \beta \in \Delta \mid \beta(Z) > 0 \} . 
\end{eqnarray}
Therefore, $\sum_{k \ge 0} \LG^k$ is a parabolic subalgebra, since 
\[
\sum_{k \ge 0} \LG^k 
= \LG_0 + \sum_{\beta(Z) \ge 0} \LG_{\beta} 
= \LG_0 + \sum_{\beta \in \Delta^+ \cup \langle \Lambda ' \rangle} \LG_{\beta} 
= \LQ_{\Lambda'} . 
\]

Conversely, let $\LQ$ be a parabolic subalgebra of $\LG$. 
By Proposition \ref{Knapp}, 
one can assume that $\LQ$ is determined by $\Lambda' \subset \Lambda$. 
Put $\Lambda \setminus \Lambda ' = \{ \alpha_{i_1}, \ldots, \alpha_{i_k} \}$ 
and define
\[
Z := H^{i_1} + \cdots + H^{i_k} . 
\]
Let $\LG = \sum \LG^k$ be the gradation with characteristic element $Z$. 
One can easily see from (\ref{eq:root1}) and (\ref{eq:root2}) that 
$\LQ = \sum_{k \ge 0} \LG^k$. 
\end{proof}

We now introduce the Langlands decomposition of a parabolic subalgebra. 

\begin{Def}
Let $\LQ_{\Lambda'}$ be the parabolic subalgebra defined by 
$\Lambda' \subset \Lambda$. 
The decomposition 
$\LQ_{\Lambda'} = \LM_{\Lambda'} + \LA_{\Lambda'} + \LN_{\Lambda'}$ 
defined by the following is called the \textit{Langlands decomposition}: 
\begin{enumerate}
\item
$\LA_{\Lambda'} 
:= 
\bigcap_{\beta \in \langle \Lambda ' \rangle} H_{\beta}^{\perp}$, 
\item
$\LM_{\Lambda'} 
:= 
(\LG_0 \ominus \LA_{\Lambda'}) + 
\sum_{\beta \in \langle \Lambda ' \rangle} \LG_{\beta}$, 
\item
$\LN_{\Lambda'} 
:= 
\sum_{\beta \in \Delta^+ \setminus \langle \Lambda ' \rangle} \LG_{\beta}$. 
\end{enumerate}
\end{Def}

The Langlands decomposition is a decomposition of a parabolic subalgebra 
into a reductive Lie algebra $\LM_{\Lambda'}$ and 
a solvable Lie algebra 
$\LS_{\Lambda '} = \LA_{\Lambda '} + \LN_{\Lambda '}$. 
We can describe Langlands decompositions in terms of gradations. 

\begin{Prop}
\label{prop:la=sum}
Let $Z := c_{i_1} H^{i_1} + \cdots + c_{i_k} H^{i_k}$ with 
$c_{i_1}, \ldots, c_{i_k} \in \Z_{>0}$. 
Denote by $\LG = \sum \LG^k$ the gradation 
with characteristic element $Z$, 
and $\LQ_{\Lambda '} = \sum_{k \ge 0} \LG^k$ the associated parabolic subalgebra, 
where 
$\Lambda ' = \{ \alpha_i \in \Lambda \mid \alpha_i(Z)=0 \}$. 
Then, for the Langlands decomposition 
$\LQ_{\Lambda'} = \LM_{\Lambda'} + \LA_{\Lambda'} + \LN_{\Lambda'}$, 
we have 
\begin{enumerate}
\item
$\LA_{\Lambda '} = \R H^{i_1} \oplus \cdots \oplus \R H^{i_k}$, 
\item
$\LM_{\Lambda '} = \LG^0 \ominus \LA_{\Lambda '}$, 
\item
$\LN_{\Lambda '} = \sum_{k>0} \LG^k$. 
\end{enumerate}
\end{Prop}

\begin{proof}
The claim (1) follows from 
$\Lambda ' = \Lambda \setminus \{ \alpha_{i_1}, \ldots, \alpha_{i_k} \}$ 
and the definition of the dual basis $\{ H^i \}$. 
The equation (\ref{eq:root1}) yields 
$\LG^0 = \LG_0 + \sum_{\beta \in \langle \Lambda ' \rangle} \LG_{\beta}$, 
which concludes (2). 
The claim (3) follows immediately from the equation (\ref{eq:root2}). 
\end{proof}

Our main theorem states that all the solvable Lie groups 
with Lie algebras $\LS_{\Lambda '} = \LA_{\Lambda '} + \LN_{\Lambda '}$ 
admit left invariant Einstein metrics. 

\section{Construction of our solvmanifolds}
\label{section_extent}

The Langlands decomposition is a decomposition of 
a parabolic subalgebra $\LQ_{\Lambda '}$ 
into reductive Lie algebra $\LM_{\Lambda '}$ and 
a solvable Lie algebra $\LS_{\Lambda '} = \LA_{\Lambda '} + \LN_{\Lambda '}$. 
In this section, we define an inner product on 
the solvable Lie algebra $\LS_{\Lambda '} = \LA_{\Lambda '} + \LN_{\Lambda '}$ 
and study the obtained solvmanifold. 

\begin{Def}
Let 
$\LQ_{\Lambda'} = \LM_{\Lambda'} + \LA_{\Lambda'} + \LN_{\Lambda'}$ 
be the Langlands decomposition of the parabolic subalgebra $\LQ_{\Lambda'}$ 
defined by $\Lambda' \subset \Lambda$. 
The solvable Lie algebra 
$\LS_{\Lambda '} = \LA_{\Lambda '} + \LN_{\Lambda '}$ 
endowed with the following inner product $\inner{}{}$ 
is called the \textit{attached solvmanifold}: 
\[
\inner{}{} := 2 B_{\sigma} |_{\LA_{\Lambda'} \times \LA_{\Lambda'}} 
+ B_{\sigma} |_{\LN_{\Lambda'} \times \LN_{\Lambda'}} . 
\]
\end{Def}

We identify the metric solvable Lie algebra 
$(\LS_{\Lambda '} = \LA_{\Lambda '} + \LN_{\Lambda '}, \inner{}{})$ 
with the simply-connected Lie group with Lie algebra $\LS_{\Lambda '}$ 
endowed with the induced left-invariant Riemannian metric. 

\begin{Prop}
The solvmanifolds 
$(\LS_{\Lambda '} = \LA_{\Lambda '} + \LN_{\Lambda '}, \inner{}{})$ 
are of Iwasawa-type. 
\end{Prop}

\begin{proof}
It is obvious that $\LA_{\Lambda '}$ is abelian. 
For every $A \in \LA_{\Lambda '}$, Proposition \ref{prop:inner} yields that 
$\ad_A$ is symmetric. 
We have 
$\LA_{\Lambda'} = \R H^{i_1} \oplus \cdots \oplus \R H^{i_k}$ 
by Proposition \ref{prop:la=sum}. 
Therefore, $A_0 := H^{i_1} + \cdots + H^{i_k}$ satisfies that 
$\ad_{A_0} |_{\LN_{\Lambda '}}$ is positive definite. 
\end{proof}

We describe the mean curvature vectors of our solvmanifolds, 
which we need for calculating the Ricci curvatures. 

\begin{Prop}
\label{mcv}
Let us consider the solvmanifold 
$(\LS_{\Lambda '} = \LA_{\Lambda '} + \LN_{\Lambda '}, \inner{}{})$, 
and define $H_0 := (1/2) \sum [\sigma(E'_j), E'_j]$, 
where $\{ E'_j \}$ is an orthonormal basis of $\LN_{\Lambda '}$. 
Then, we have 
\begin{enumerate}
\item
$H_0 \in \LA_{\Lambda '}$, and 
\item
the mean curvature vector of $(\LS_{\Lambda '}, \inner{}{})$ 
coincides with $H_0$. 
\end{enumerate}
\end{Prop}

\begin{proof}
We can and do assume that each $E'_j$ is contained in a root space. 
Since $[\sigma(E'_j), E'_j]$ is a root vector, we have 
\[
2 H_0 = \sum_{\alpha(Z)>0} (\dim \LG_{\alpha}) H_{\alpha} , 
\]
where $Z$ is the characteristic element. 
This concludes $H_0 \in \LA$. 
One can represent $Z = c_{i_1} H^{i_1} + \cdots + c_{i_k} H^{i_k}$ with 
$c_{i_1}, \ldots, c_{i_k} \in \Z_{>0}$. 
It follows from Proposition \ref{prop:la=sum} that 
\[
\LA \ominus \LA_{\Lambda '} 
= \mathrm{span} \{ H_{\alpha_j} \mid j \neq i_1, \ldots, i_k \} . 
\]
Therefore, in order to prove (1), 
we have only to show $B_{\sigma}(H_{\alpha_j}, H_0) = 0$ 
for every $j \neq i_1, \ldots, i_s$. 
Let us take $j \neq i_1, \ldots, i_s$. 
Let $s_j$ be the reflection with respect to $\alpha_j^{\perp}$, 
which is an element of the Weyl group. 
One can see that 
\begin{eqnarray*}
s_j(H_0) 
= s_j(\sum_{\alpha(Z)>0} (\dim \LG_{\alpha}) H_{\alpha}) 
= \sum_{\alpha(Z)>0} (\dim \LG_{\alpha}) H_{s_j(\alpha)} . 
\end{eqnarray*}
One knows $\dim \LG_{\alpha} = \dim \LG_{s_j(\alpha)}$, 
since the multiplicities of roots are invariant 
under the action of the Weyl group. 
Furthermore, 
since $s_j$ acts trivially on $\LA_{\Lambda '}$, 
one has $s_j(Z)=Z$, and hence $s_j$ preserves 
$\{ \alpha \in \Delta \mid \alpha(Z)>0 \}$. 
Therefore we conclude that 
\begin{eqnarray*}
s_j(H_0) 
= \sum_{\alpha(Z)>0} (\dim \LG_{s_j(\alpha)}) H_{s_j(\alpha)} 
= H_0 . 
\end{eqnarray*}
This means $B_{\sigma}(H_{\alpha_j}, H_0) = 0$, 
which proves (1). 
In order to show (2), let $A \in \LA_{\Lambda '}$. 
Since $\LA_{\Lambda '}$ is abelian, one can see that 
\begin{eqnarray*}
\tr(\ad_A) 
= \sum \inner{\ad_A(E'_j)}{E'_j} 
= \sum B_{\sigma}([A, E'_j], E'_j) . 
\end{eqnarray*}
It follows from Proposition \ref{prop:inner} and the definition of our inner product 
that 
\begin{eqnarray*}
\sum B_{\sigma}([A, E'_j], E'_j) 
= \sum B_{\sigma}(A, [\sigma(E'_j), E'_j]) 
= B_{\sigma}(A, 2 H_0)
= \inner{A}{H_0} . 
\end{eqnarray*}
Therefore, $H_0$ is the mean curvature vector. 
\end{proof}

This concludes that, the rank one reductions 
$(\R H_0 + \LN_{\Lambda '}, \inner{}{})$ 
of our solvmanifolds are nothing but the solvmanifolds studied in \cite{T7}. 
Note that, it is not easy to give the mean curvature vectors explicitly. 

\medskip 

Our next aim is to see that, 
for the minimal parabolic subalgebra $\LQ_{\emptyset}$, 
the attached solvmanifold 
$(\LS_{\emptyset}, \inner{}{})$ 
coincides with the symmetric space associated with the pair $(\LG, \LK)$. 

\begin{Prop}
\label{prop:symm}
For the minimal parabolic subalgebra $\LQ_{\emptyset}$ of $\LG$, 
the attached solvmanifold 
$(\LS_{\emptyset} = \LA_{\emptyset} + \LN_{\emptyset}, \inner{}{})$ 
is a symmetric space and an Einstein manifold with Einstein constant $-(1/4)$. 
\end{Prop}

\begin{proof}
Let $\sigma$ be the Cartan involution of $\LG$ 
and $\LG = \LK + \LP$ the Cartan decomposition. 
The pair $(\LG, \LK)$ gives the 
symmetric space of noncompact type $(G/K, g)$, 
where the metric $g$ is defined by $B$. 
It is known that $(G/K, g)$ 
is an Einstein manifold with Einstein constant $-(1/2)$ 
(see \cite[Theorem 7.73]{Besse}). 
We will show that our solvmanifold is isometric to $(G/K, 2g)$, 
which is also a symmetric space and 
an Einstein manifold with Einstein constant $-(1/4)$. 

Let $\Lambda$ be the set of simple roots, and 
$\LG = \LK + \LA + \LN$ the Iwasawa decomposition determined by $\Lambda$. 
Recall that 
$\LA + \LN = \LA_{\emptyset} + \LN_{\emptyset}$. 
Denote by $G = KAN$ the corresponding decomposition of the group. 
Thus, one has an isomorphism 
\[
\varphi : AN \rightarrow G/K : x \mapsto [x] = x.o , 
\]
where $o$ denotes the origin. 
By identifying $\LP = T_o (G/K)$, 
the differential of $\varphi$ at the identity $e$ is given by
\[
(d \varphi)_e : \LA + \LN \rightarrow \LP : 
V \mapsto \frac{d}{dt} \exp (tV).o |_{t=0} = V_{\LP} , 
\]
where $V_{\LP}$ denotes the $\LP$-component of $V$. 
Thus, for $A \in \LA$ and $X \in \LN$, we have 
\[
(d \varphi)_e (A + X) = A + (1/2)(X - \sigma(X)) . 
\]
Note that $B(X,X')=0$. 
Therefore, we have 
\begin{eqnarray*}
&&
2 B( (d \varphi)_e (A + X), (d \varphi)_e (A' + X') ) \\ 
& = & 
2 B( A + (1/2)(X - \sigma(X)) , A' + (1/2)(X' - \sigma(X')) ) \\ 
& = & 
2 B(A, A') + (1/2) B(X, -\sigma(X')) + (1/2) B(-\sigma(X), X') \\ 
& = & 
2 B_{\sigma} (A, A') + B_{\sigma}(X, X') \\ 
& = & 
\inner{A+X}{A'+X'} . 
\end{eqnarray*}
Since the metric $2g$ is defined by $2B$, 
we conclude that $\varphi$ is an isometry. 
\end{proof}

Therefore, the class of our solvmanifolds, attached to parabolic subalgebras, 
contains all symmetric spaces of noncompact type. 
Our solvmanifolds contain many non-symmetric ones, 
but they still closely relate to symmetric spaces. 

\begin{Prop}
All of our solvmanifolds 
$(\LS_{\Lambda'} = \LA_{\Lambda '} + \LN_{\Lambda '}, \inner{}{})$ 
are naturally Riemannian submanifold of the symmetric spaces 
$(\LS_{\emptyset} = \LA_{\emptyset} + \LN_{\emptyset}, \inner{}{})$. 
\end{Prop}

The proof is easy from the construction. 
This observation is the key for the proof of our main theorem. 

\section{Proof of the main theorem}
\label{section_Einstein}

In this section, we show our main theorem 
which states that all of our solvmanifolds 
$(\LS_{\Lambda '} = \LA_{\Lambda '} + \LN_{\Lambda '}, \inner{}{})$ 
are Einstein. 
The strategy for the proof is the use of submanifold theory. 
We compare the Ricci curvatures of 
$(\LS_{\Lambda '} = \LA_{\Lambda '} + \LN_{\Lambda '}, \inner{}{})$ 
and the Ricci curvatures of the symmetric spaces 
$(\LS_{\emptyset} = \LA_{\emptyset} + \LN_{\emptyset}, \inner{}{})$. 
Denote by $\ric^{\LG}$ the Ricci curvature 
of a metric Lie algebra $(\LG, \inner{}{})$. 

\medskip 

First of all, the Ricci curvature for $\LA_{\Lambda '}$-direction 
can be calculated directly. 

\begin{Lem}
\label{lem:a}
$\ric^{\LS_{\Lambda'}}(A,A') = -(1/4) \inner{A}{A'}$ 
for every $A, A' \in \LA_{\Lambda'}$. 
\end{Lem}

\begin{proof}
Let $A, A' \in \LA_{\Lambda'}$. 
The definition of our inner product yields that 
\begin{eqnarray*}
\inner{A}{A'} 
= 2 B_{\sigma}(A,A') 
= 2 B(A,A') 
= 2 \tr (\ad^{\LG}_A) \circ (\ad^{\LG}_{A'}) . 
\end{eqnarray*}
We have an orthogonal decomposition 
$\LG = \sigma(\LN_{\Lambda '}) + \LG^0 + \LN_{\Lambda '}$, 
which is invariant under $\ad_A$ and $\ad_{A'}$. 
Therefore, the traces can also be decomposed into these three directions. 
First of all, one has 
\begin{eqnarray*}
\tr (\ad^{\LG}_A) \circ (\ad^{\LG}_{A'}) |_{\LG^0} = 0 , 
\end{eqnarray*}
since $[\LA_{\Lambda '}, \LG^0] =0$ 
(see Proposition \ref{prop:la=sum}). 
Let $\{ E'_j \}$ be an orthonormal basis of $\LN_{\Lambda '}$. 
Then, since $\LA_{\Lambda '}$ is abelian, one has 
\begin{eqnarray*}
\tr (\ad^{\LG}_A) \circ (\ad^{\LG}_{A'}) |_{\LN_{\Lambda '}}
= \sum B_{\sigma}([A, [A', E'_j], E'_j) 
= \tr (\ad_A) \circ (\ad_{A'}) . 
\end{eqnarray*}
Finally, since $\{ \sigma(E'_j) \}$ is an orthonormal basis of 
$\sigma(\LN_{\Lambda '})$, one can see that 
\begin{eqnarray*}
\tr (\ad^{\LG}_A) \circ (\ad^{\LG}_{A'}) |_{\sigma(\LN_{\Lambda '})}
= \sum B_{\sigma}([A, [A', \sigma(E'_j)]], \sigma(E'_j)) 
= \tr (\ad_{A}) \circ (\ad_{A'}) . 
\end{eqnarray*}
Altogether we obtain 
\begin{eqnarray*}
\inner{A}{A'} = 4 \tr (\ad_A) \circ (\ad_{A'}) . 
\end{eqnarray*}
Theorem \ref{wolter}, 
which states 
$\ric^{\LS_{\Lambda'}}(A,A') = - \tr (\ad_A) \circ (\ad_{A'})$, 
completes the proof. 
\end{proof}

In order to prove the main theorem, 
we have to know the Ricci curvature for $\LN_{\Lambda'}$-direction. 
We use the $(1,1)$-tensor $\Ric^{\LG}$ defined by 
\[
\ric^{\LG}(X,Y) = \inner{\Ric^{\LG}(X)}{Y} . 
\]

\begin{Lem}
\label{lem:n}
Let $\{ E_j^{\perp} \}$ be an orthonormal basis of 
$\LN_{\emptyset} \ominus \LN_{\Lambda'}$ and 
define $H_0^{\perp} := (1/2) \sum [\sigma E_j^{\perp}, E_j^{\perp}]$. 
Then, for every $X \in \LN_{\Lambda'}$, we have 
\[
\Ric^{\LN_{\emptyset}}(X) - \Ric^{\LN_{\Lambda'}}(X) = [H_0^{\perp}, X] . 
\] 
\end{Lem}

\begin{proof}
For $X, Y \in \LN_{\Lambda '}$, 
Proposition \ref{prop:inner} shows that 
$(\ad_{X})^{\ast} Y = -[\sigma X, Y]_{\LN_{\Lambda '}}$, 
where subscript $\LN_{\Lambda '}$ denotes the $\LN_{\Lambda '}$-component. 
Hence, Proposition \ref{ricci-nil} yields that 
\begin{eqnarray*}
\Ric^{\LN_{\emptyset}} (X) & = & 
- (1/4) \sum [E_i, [\sigma E_i, X]_{\LN_{\emptyset}}] 
+ (1/2) \sum [\sigma E_i, [E_i, X]]_{\LN_{\emptyset}}, \\
\Ric^{\LN_{\Lambda'}} (X) & = & 
- (1/4) \sum [E'_j, [\sigma E'_j, X]_{\LN_{\Lambda'}}] 
+ (1/2) \sum [\sigma E'_j, [E'_j, X]]_{\LN_{\Lambda'}}, 
\end{eqnarray*}
where $\{ E'_j \}$ is an orthonormal basis of $\LN_{\Lambda'}$ 
and $\{ E_i \} := \{ E'_j \} \cup \{ E_k^{\perp} \}$. 
We can and do assume that each $E_i$ is contained in a root space. 
For simplicity of the notations, put 
\begin{eqnarray*}
& A := \Sum [E_k^{\perp}, [\sigma E_k^{\perp}, X]_{\LN_{\emptyset}}] , 
& B := \sum [E'_j, [\sigma E'_j, X]_{\LN_{\emptyset} \ominus \LN_{\Lambda '}}] \\ 
& C := \Sum [\sigma E_k^{\perp}, [E_k^{\perp}, X]]_{\LN_{\emptyset}} , 
& D := \sum [\sigma E'_j, [E'_j, X]]_{\LN_{\emptyset} \ominus \LN_{\Lambda '}} . 
\end{eqnarray*}
Therefore, it is easy to see that 
\begin{eqnarray*}
\sum [E_i, [\sigma E_i, X]_{\LN_{\emptyset}}] 
& = & 
A + B + \sum [E'_j, [\sigma E'_j, X]_{\LN_{\Lambda'}}] , \\ 
\sum [\sigma E_i, [E_i, X]]_{\LN_{\emptyset}}
& = & 
C + D + \sum [\sigma E'_j, [E'_j, X]]_{\LN_{\Lambda'}} . 
\end{eqnarray*}
We then obtain 
\begin{eqnarray}
\label{eq:diff}
\Ric^{\LN_{\emptyset}} (X) - \Ric^{\LN_{\Lambda'}} (X) 
= - (1/4) A - (1/4) B + (1/2) C + (1/2) D . 
\end{eqnarray}

Our first claim is 
\begin{eqnarray}
\label{eq:claim}
[\sigma E_k^{\perp}, X] \in \LN_{\Lambda '} . 
\end{eqnarray}
Denote by $Z$ the characteristic element. 
Recall that $\LN_{\Lambda '} = \sum_{\alpha(Z)>0} \LG_{\alpha}$ 
and hence $\LN_{\emptyset} \ominus \LN_{\Lambda '} = 
\sum_{\beta > 0, \beta(Z)=0} \LG_{\beta}$. 
Therefore we have $E_k^{\perp} \in \LG_{\beta}$ with $\beta >0$ and $\beta(Z)=0$. 
From the linearity, we have only to consider the case 
$X \in \LG_{\alpha}$ with $\alpha(Z) > 0$. 
One can see that 
$[\sigma E_k^{\perp}, X] \in \LG_{\alpha - \beta}$ and $(\alpha - \beta)(Z) >0$. 
This concludes 
$[\sigma E_k^{\perp}, X] \in \LN_{\Lambda '}$. 

The claim (\ref{eq:claim}) immediately yields that 
\begin{eqnarray}
\label{eq:a}
A = \sum [E_k^{\perp}, [\sigma E_k^{\perp}, X]] . 
\end{eqnarray}

Our second claim is 
\begin{eqnarray}
\label{eq:b}
A = B . 
\end{eqnarray}
Since $\{ E_k^{\perp} \}$ is an orthonormal basis of 
$\LN_{\emptyset} \ominus \LN_{\Lambda '}$, one has 
\begin{eqnarray*}
B = \sum_j [E'_j, \sum_k \inner{[\sigma E'_j, X]}{E_k^{\perp}} E_k^{\perp}] . 
\end{eqnarray*}
By the property of our inner product, we have 
\begin{eqnarray*}
\inner{[\sigma E'_j, X]}{E_k^{\perp}} 
= - \inner{X}{[E'_j, E_k^{\perp}]} 
= \inner{[X, \sigma E_k^{\perp}]}{E'_j} . 
\end{eqnarray*}
Since $\{ E'_j \}$ is an orthonormal basis of $\LN_{\Lambda '}$, 
one can see that 
\begin{eqnarray*}
B = \sum_{j, k} [ \inner{[X, \sigma E_k^{\perp}]}{E'_j} E'_j, E_k^{\perp}] 
= \sum_{k} [ [X, \sigma E_k^{\perp}]_{\LN_{\Lambda '}}, E_k^{\perp}] . 
\end{eqnarray*}
Therefore we conclude $B=A$ from (\ref{eq:claim}) and (\ref{eq:a}). 

Next, we calculate $C$. 
A similar argument to the proof for (\ref{eq:claim}) implies that 
\[ 
C = \sum [\sigma E_k^{\perp}, [E_k^{\perp}, X]] . 
\]
Thus, the Jacobi identity yields that 
\begin{eqnarray}
\label{eq:c}
C = - \sum [E_k^{\perp}, [X, \sigma E_k^{\perp}]] 
- \sum [X, [ \sigma E_k^{\perp}, E_k^{\perp}]] 
= A + [2 H_0^{\perp}, X] . 
\end{eqnarray}

Our last claim is 
\begin{eqnarray}
\label{eq:d}
D = 0 . 
\end{eqnarray}
By assumption, we have $E'_j \in \LG_{\alpha}$ with $\alpha(Z)>0$. 
Furthermore, from the linearity, we have only to consider the case 
$X \in \LG_{\gamma}$ with $\gamma(Z)>0$. 
Then, it is obvious that 
\[
[\sigma E'_j, [E'_j, X]] \in \LG_{\gamma} \subset \LN_{\Lambda '} , 
\]
which can not have $\LN_{\emptyset} \ominus \LN_{\Lambda '}$-component. 
This concludes the claim (\ref{eq:d}). 

We conclude the lemma 
by combining (\ref{eq:diff}), (\ref{eq:a}), (\ref{eq:b}), (\ref{eq:c}) 
and (\ref{eq:d}). 
\end{proof}

Lemma \ref{lem:n} gives the relations between 
$\ric^{\LN_{\emptyset}}$ and $\ric^{\LN_{\Lambda '}}$, 
from which we can obtain the relations between 
$\ric^{\LS_{\emptyset}}$ and $\ric^{\LS_{\Lambda '}}$. 

\begin{Thm}
\label{thm:main}
Let $\LG$ be a semisimple Lie algebra, 
$\LQ_{\Lambda '}$ a parabolic subalgebra of $\LG$, and 
$(\LS_{\Lambda '} = \LA_{\Lambda '} + \LN_{\Lambda '}, \inner{}{})$ 
the attached solvmanifold. 
Then, the Ricci curvatures satisfy 
$\ric^{\LS_{\Lambda '}} = \ric^{\LS_{\emptyset}}$ on 
$\LS_{\Lambda '} \times \LS_{\Lambda '}$. 
Therefore, every 
$(\LS_{\Lambda '} = \LA_{\Lambda '} + \LN_{\Lambda '}, \inner{}{})$ 
is an Einstein manifold with Einstein constant $-(1/4)$. 
\end{Thm}

\begin{proof}
Let $A, A' \in \LA_{\Lambda '}$. 
In this case, Lemma \ref{lem:a} yields that 
\[
\ric^{\LS_{\Lambda '}}(A,A') 
= -(1/4) \inner{A}{A'} 
= \ric^{\LS_{\emptyset}}(A,A') . 
\]
Let $A \in \LA_{\Lambda '}$ and $X \in \LN_{\Lambda '}$. 
It follows directly from Theorem \ref{wolter} that 
\[
\ric^{\LS_{\Lambda '}}(A,X) 
= 0
= \ric^{\LS_{\emptyset}}(A,X) . 
\]
Finally, let us consider the case $X, Y \in \LN_{\Lambda '}$. 
As in the proof for Lemma \ref{lem:n}, 
let $\{ E'_j \}$ and $\{ E^{\perp}_k \}$ be orthonormal basis of 
$\LN_{\Lambda '}$ and $\LN_{\emptyset} \ominus \LN_{\Lambda '}$, respectively. 
Define $H'_0 = (1/2) \sum [\sigma(E'_j), E'_j]$, 
$H^{\perp}_0 = (1/2) \sum [\sigma(E^{\perp}_k), E^{\perp}_k]$, 
and $H_0 = H'_0 + H^{\perp}_0$. 
Therefore, Proposition \ref{mcv} yields that 
$H'_0$ and $H_0$ are the mean curvature vectors of 
$(\LS_{\Lambda '}, \inner{}{})$ and $(\LS_{\emptyset}, \inner{}{})$, respectively. 
By Theorem \ref{wolter}, one has 
\begin{eqnarray*}
\ric^{\LS_{\emptyset}}(X,Y) 
& = & 
\inner{\Ric^{\LN_{\emptyset}}(X)}{Y} - \inner{[H_0, X]}{Y} , \\ 
\ric^{\LS_{\Lambda '}}(X,Y) 
& = & 
\inner{\Ric^{\LN_{\Lambda '}}(X)}{Y} - \inner{[H'_0, X]}{Y} .
\end{eqnarray*}
From Lemma \ref{lem:n}, we have 
\begin{eqnarray*}
&&
\ric^{\LS_{\emptyset}}(X,Y) - \ric^{\LS_{\Lambda '}}(X,Y) \\
& = & 
\inner{\Ric^{\LN_{\emptyset}}(X) - \Ric^{\LN_{\Lambda '}}(X)}{Y} 
- \inner{[H_0 - H'_0, X]}{Y} \\ 
& = & 
\inner{[H_0^{\perp}, X]}{Y} 
- \inner{[H^{\perp}_0, X]}{Y} \\ 
& = & 0 . 
\end{eqnarray*}
This completes the proof. 
The Einstein constant comes from Proposition \ref{prop:symm}. 
\end{proof}

Our Einstein solvmanifolds are Riemannian submanifolds of 
symmetric spaces of noncompact type, 
and the Ricci curvatures coincide with 
the restrictions of the Ricci curvatures of the ambient spaces. 
This would be a remarkable property of Riemannian submanifolds. 
The rank one reductions of standard Einstein solvmanifolds 
$(\LS = \LA + \LN, \inner{}{})$ 
are also examples of submanifolds with this property. 
In this case, 
one can remove a particular subspace in $\LA$, keeping Ricci curvature. 
Our solvmanifolds state that, 
one can also remove a particular subspace in $\LS$, 
not necessarily contained in $\LA$, keeping Ricci curvature. 
The study of submanifolds with this property seems to be interesting, 
which may provide new examples of Einstein solvmanifolds.

\section{Minimality}
\label{sec:minimal}

In this section, we show that all of our solvmanifolds 
$(\LS_{\Lambda '} = \LA_{\Lambda '} + \LN_{\Lambda '}, \inner{}{})$ 
are minimal, but not totally geodesic, homogeneous submanifolds 
of the symmetric spaces of noncompact type 
$(\LS_{\emptyset} = \LA_{\emptyset} + \LN_{\emptyset}, \inner{}{})$. 
We need the second fundamental form 
$h : \LS_{\Lambda '} \times \LS_{\Lambda '} \rightarrow 
\LS_{\emptyset} \ominus \LS_{\Lambda '}$, 
defined by 
\[
h(X,Y) = \nabla^{\LS_{\emptyset}}_X Y - \nabla^{\LS_{\Lambda '}}_X Y , 
\]
where 
$\nabla^{\LS_{\emptyset}}$ and $\nabla^{\LS_{\Lambda '}}$ 
denote the Levi-Civita connections of 
$(\LS_{\Lambda '}, \inner{}{})$ and $(\LS_{\emptyset}, \inner{}{})$, 
respectively. 

\begin{Lem}
\label{lem:h}
Let $h$ be the second fundamental of the submanifold 
$(\LS_{\Lambda '} = \LA_{\Lambda '} + \LN_{\Lambda '}, \inner{}{})$ 
in 
$(\LS_{\emptyset} = \LA_{\emptyset} + \LN_{\emptyset}, \inner{}{})$. 
Then we have 
\begin{enumerate}
\item
$h(A,A)=0$ for $A \in \LA_{\Lambda '}$, 
\item
$h(X_{\alpha}, X_{\alpha}) = 
(1/2) [\sigma(X_{\alpha}), X_{\alpha}]_{\LA_{\emptyset} \ominus \LA_{\Lambda '}}$ 
for $X_{\alpha} \in \LG_{\alpha} \subset \LN_{\Lambda '}$. 
\end{enumerate}
\end{Lem}

\begin{proof}
First of all, one has 
\[
h(X,Y) = U^{\LS_{\emptyset}}(X,Y) - U^{\LS_{\Lambda '}}(X,Y) . 
\]
Recall that 
$U^{\LS_{\Lambda '}} : \LS_{\Lambda '} \times \LS_{\Lambda '} 
\rightarrow \LS_{\Lambda '}$ 
is defined by 
\[
2 \inner{U^{\LS_{\Lambda '}}(X,Y)}{Z} 
= \inner{[Z,X]}{Y} + \inner{X}{[Z,Y]} . 
\]
Let $A \in \LA_{\Lambda '}$. 
Thus, it is easy to see that 
$U^{\LS_{\emptyset}}(A,A) = 0$ and $U^{\LS_{\Lambda '}}(A,A) = 0$. 
This proves (1). 
In order to prove (2), 
let us take $X_{\alpha} \in \LG_{\alpha} \subset \LN_{\Lambda '}$. 
For every $Z \in \LS_{\Lambda '}$, one has 
\begin{eqnarray*}
\inner{U^{\LS_{\Lambda '}}(X_{\alpha},X_{\alpha})}{Z} 
& = & 
\inner{[Z,X_{\alpha}]}{X_{\alpha}} \\
& = & 
B_{\sigma}([Z,X_{\alpha}], X_{\alpha}) \\
& = & 
B_{\sigma}(Z, [\sigma(X_{\alpha}),X_{\alpha}]) \\
& = & 
(1/2) \inner{Z}{[\sigma(X_{\alpha}),X_{\alpha}]} . 
\end{eqnarray*}
Note that $[\sigma(X_{\alpha}),X_{\alpha}] = |X_{\alpha}|^2 H_{\alpha} \in \LA$. 
Thus, we obtain 
\begin{eqnarray*}
U^{\LS_{\emptyset}}(X_{\alpha},X_{\alpha}) 
& = & 
(1/2) [\sigma(X_{\alpha}),X_{\alpha}]_{\LA_{\emptyset}} , \\ 
U^{\LS_{\Lambda '}}(X_{\alpha},X_{\alpha}) 
& = & 
(1/2) [\sigma(X_{\alpha}),X_{\alpha}]_{\LA_{\Lambda '}} , 
\end{eqnarray*}
which conclude (2). 
\end{proof}

Now we show the minimality, that is, the trace of $h$ vanishes. 

\begin{Thm}
All of our solvmanifolds 
$(\LS_{\Lambda '} = \LA_{\Lambda '} + \LN_{\Lambda '}, \inner{}{})$ 
are minimal submanifolds in 
$(\LS_{\emptyset} = \LA_{\emptyset} + \LN_{\emptyset}, \inner{}{})$. 
\end{Thm}

\begin{proof}
Let $\{ A_i \}$ and $\{ E_j \}$ are orthonormal basis of 
$\LA_{\Lambda '}$ and $\LN_{\Lambda '}$, respectively. 
We may and do assume that each $E_j$ is contained in a root space. 
Denote by $H$ the mean curvature vector of our submanifold. 
Thus, Lemma \ref{lem:h} yields that 
\begin{eqnarray*}
H = \sum h(A_i, A_i) + \sum h(E_j, E_j) 
= (1/2) \sum [\sigma(E_j), E_j]_{\LA_{\emptyset} \ominus \LA_{\Lambda '}} . 
\end{eqnarray*}
We have proved in Proposition \ref{mcv} that 
\[
(1/2) \sum [\sigma(E_j), E_j] \in \LA_{\Lambda '} . 
\]
Therefore, we conclude $H=0$. 
\end{proof}

Lemma \ref{lem:h} also implies that our submanifolds are not totally geodesic, 
that is, $h$ do not vanish. 
In order to show this, we have to exclude trivial cases. 

\begin{Def}
A subset $\Lambda ' \subset \Lambda$ is said to be \textit{trivial} 
if $\Lambda '$ and $\Lambda \setminus \Lambda '$ are orthogonal. 
\end{Def}

It is obvious that $\Lambda ' = \emptyset$ is trivial. 
Assume that $\Lambda '$ is trivial. 
Thus, the root system has an orthogonal decomposition 
$\Delta = \langle \Lambda ' \rangle \sqcup 
\langle \Lambda \setminus \Lambda ' \rangle$, 
and hence we have an orthogonal decomposition into corresponding ideals, 
$\LG = \LG_{\langle \Lambda ' \rangle} 
+ \LG_{\langle \Lambda \setminus \Lambda ' \rangle}$. 
Therefore, the attached solvmanifold 
$(\LS_{\Lambda '}, \inner{}{})$ 
coincides with the symmetric space associated with
$\LG_{\langle \Lambda \setminus \Lambda ' \rangle}$, 
which is totally geodesic. 

\begin{Prop}
\label{prop:nontotallygeodesic}
Assume that $\Lambda '$ is not trivial. 
Then, the solvmanifolds 
$(\LS_{\Lambda '} = \LA_{\Lambda '} + \LN_{\Lambda '}, \inner{}{})$ 
are not totally geodesic in 
$(\LS_{\emptyset} = \LA_{\emptyset} + \LN_{\emptyset}, \inner{}{})$. 
\end{Prop}

\begin{proof}
By assumption, there exist 
$\alpha \in \Lambda \setminus \Lambda '$ and $\beta \in \Lambda '$ 
such that $\inner{\alpha}{\beta} \neq 0$. 
We have 
$\LG_{\alpha} \in \LN_{\Lambda '}$ and 
$\LG_{\beta} \subset \LN_{\emptyset} \ominus \LN_{\Lambda '}$. 
Let $X_{\alpha}$ be a unit vector in $\LG_{\alpha}$. 
Then, Lemma \ref{lem:h} yields that 
\begin{eqnarray*}
2 h(X_{\alpha},X_{\alpha}) 
= [\sigma(X_{\alpha}), X_{\alpha}]_{\LA_{\emptyset} \ominus \LA_{\Lambda '}} 
= (H_{\alpha})_{\LA_{\emptyset} \ominus \LA_{\Lambda '}} . 
\end{eqnarray*}
Therefore, we obtain $h(X_{\alpha}, X_{\alpha}) \neq 0$, 
because 
$H_{\beta} \in \LA_{\emptyset} \ominus \LA_{\Lambda '}$ and 
$\inner{H_{\alpha}}{H_{\beta}} = \inner{\alpha}{\beta} \neq 0$. 
\end{proof}

\bibliographystyle{amsplain}

\begin{thebibliography}{99} 

\bibitem{A} 
D.\ V.\ Alekseevskii, 
Homogeneous Riemannian spaces of negative curvature. 
\textit{Mat. Sb.} {\bf 25} (1975), 87-109; 
English translation, {\it Math. USSR-Sb.} {\bf 96} (1975), 93--117. 

\bibitem{BTV} 
J.\ Berndt, F.\ Tricerri, L.\ Vanhecke, 
\textit{Generalized Heisenberg groups and Damek-Ricci harmonic spaces}. 
Lect. Notes in Math. {\bf 1598} (1995), 
Springer-Verlag Berlin Heidelberg. 

\bibitem{Besse} 
A.\ Besse, 
\textit{Einstein manifolds}. 
Ergeb. Math. {\bf 10} (1987), Springer-Verlag, Berlin-Heidelberg. 

\bibitem{DR} 
E.\ Damek, F.\ Ricci, 
Harmonic analysis on solvable extensions of $H$-type groups.
\textit{J. Geom. Anal.} {\bf 2} (1992), 213--248. 

\bibitem{Eb} 
P.\ B.\ Eberlein, 
\textit{Geometry of nonpositively curved manifolds}. 
University of Chicago Press, Chicago, London, 1996. 


\bibitem{H} 
J.\ Heber, 
Noncompact homogeneous Einstein spaces. 
\textit{Invent. Math.} {\bf 133} (1998), 279-352. 


\bibitem{Hel} 
S.\ Helgason, 
\textit{Differential geometry, Lie groups, and symmetric spaces}. 
Graduate Studies in Mathematics {\bf 34}, 
American Mathematical Society, Providence, RI, 2001. 

\bibitem{KA} 
S.\ Kaneyuki, H.\ Asano, 
Graded Lie algebras and generalized Jordan triple systems. 
\textit{Nagoya Math. J.} {\bf 112} (1988), 81--115. 

\bibitem{Knapp}
A.\ W.\ Knapp, 
\textit{Lie groups beyond an introduction}. 
Second edition. Progress in Mathematics, 140. 
Birkhauser Boston, Inc., Boston, MA, 2002. xviii+812 pp. 

\bibitem{L-dga} 
J.\ Lauret, 
Minimal metrics on nilmanifolds. 
In: Bure\v s, J., Kowalski, O., Krupka, D., Slovak, J. (eds.)
Differential geometry and its applications
Proceedings of the 9th International Conference (DGA 2004) held in Prague, 
August 30--September 3, 2004, 
79--97. 
Matfyzpress, Prague (2005). 

\bibitem{L-standard} 
J.\ Lauret, 
Einstein solvmanifolds are standard. 
Preprint (2007). 
\textit{arXiv}:math/0703472. 

\bibitem{Mori} 
K.\ Mori, 
Einstein metrics on Boggino-Damek-Ricci-type solvable Lie groups. 
\textit{Osaka J. Math.} {\bf 39} (2002), 345--362. 

\bibitem{Niko} 
Y.\ Nikolayevsky, 
Einstein solvmanifolds with free nilradical. 
\textit{Ann. Glob. Anal. Geom.}, to appear. 
\textit{arXiv}:math/0610020. 

\bibitem{OV} 
A.\ L.\ Onishchik and E.\ B.\ Vinberg, 
\textit{Lie Groups and Lie Algebras III}. 
Springer-Verlag. 

\bibitem{T4} 
H.\ Tamaru, 
On certain subalgebras of graded Lie algebras. 
\textit{Yokohama Math. J.} {\bf 46} (1999), 127--138. 

\bibitem{T7} 
H.\ Tamaru, 
Noncompact homogeneous Einstein manifolds attached to graded Lie algebras. 
\textit{Math. Z.}, to appear. 
\textit{arXiv}: math.DG/0610675. 

\bibitem{T-dga} 
H.\ Tamaru, 
A class of noncompact homogeneous Einstein manifolds. 
In: Bure\v s, J., Kowalski, O., Krupka, D., Slovak, J. (eds.)
Differential geometry and its applications, 
Proceedings of the 9th International Conference (DGA 2004) held in Prague, 
August 30--September 3, 2004, 119--127. 
Matfyzpress, Prague (2005). 

\bibitem{W} 
T.\ H.\ Wolter, 
Einstein metrics on solvable Lie groups. 
\textit{Math. Z.} \textbf{206} (1991), 457--471. 

\end{thebibliography}

\end{document}